\documentclass{article}
\usepackage{stmaryrd}
\usepackage{mathrsfs}
\usepackage[centertags]{amsmath}
\usepackage{amsfonts, dsfont}
\usepackage{amssymb}
\usepackage{amsthm}
\usepackage{color}

\input xy
\xyoption{all}

\theoremstyle{plain}
\newtheorem{thm}{Theorem}[section]
\newtheorem{cor}[thm]{Corollary}
\newtheorem{lem}[thm]{Lemma}
\newtheorem{prop}[thm]{Proposition}
\newtheorem*{thm2}{Theorem}

\theoremstyle{definition}
\newtheorem{defn}[thm]{Definition}
\newtheorem{remark}[thm]{Remark}
\newtheorem*{ack}{Acknowledgments}

\newcommand{\bd}{\begin{defn}}
\newcommand{\ed}{\end{defn}}
\newcommand{\bl}{\begin{lem}}
\newcommand{\el}{\end{lem}}
\newcommand{\bp}{\begin{prop}}
\newcommand{\ep}{\end{prop}}
\newcommand{\bt}{\begin{thm}}
\newcommand{\et}{\end{thm}}
\newcommand{\bc}{\begin{cor}}
\newcommand{\ec}{\end{cor}}
\newcommand{\br}{\begin{remark}}
\newcommand{\er}{\end{remark}}
\newcommand{\bdi}{\begin{diagram}}
\newcommand{\edi}{\end{diagram}}
\newcommand{\beq}{\begin{equation}}
\newcommand{\eeq}{\end{equation}}
\newcommand{\ba}{\begin{array}}
\newcommand{\ea}{\end{array}}
\newcommand{\bpf}{\begin{proof}}
\newcommand{\epf}{\end{proof}}

\newcommand{\Z}{\mathds{Z}}
\newcommand{\Q}{\mathds{Q}}
\newcommand{\Zp}{\mathds{Z}_{p}}

\newcommand{\al}{\alpha}

\newcommand{\Ga}{\Gamma}

\newcommand{\la}{\lambda}

\newcommand{\Cl}{\mathrm{Cl}}

 \DeclareMathOperator{\Gal}{Gal}
 \DeclareMathOperator{\rank}{rank}

\DeclareMathOperator{\Ext}{Ext}

\newcommand{\X}{\mathcal{X}}
\newcommand{\Y}{\mathcal{Y}}

\newcommand{\M}{\mathfrak{M}}

\newcommand{\ot}{\otimes}

\newcommand{\plim}{\displaystyle \mathop{\varprojlim}\limits}

\newcommand{\lra}{\longrightarrow}

\newcommand{\ps}[1]{\llbracket #1 \rrbracket}

\begin{document}
\title{A note on asymptotic class number upper bounds in $p$-adic Lie extensions}
 \author{
  Meng Fai Lim\footnote{School of Mathematics and Statistics $\&$ Hubei Key Laboratory of Mathematical Sciences,
Central China Normal University, Wuhan, 430079, P.R.China.
 E-mail: \texttt{limmf@mail.ccnu.edu.cn}} }
\date{}
\maketitle

\begin{abstract} \footnotesize
\noindent Let $p$ be an odd prime and $F_{\infty}$ a $p$-adic Lie extension of a number field $F$ with Galois group $G$. Suppose that $G$ is a compact pro-$p$ $p$-adic Lie group with no torsion and that it contains a closed normal subgroup $H$ such that $G/H\cong \Zp$. Under various assumptions, we establish asymptotic upper bounds for the growth of $p$-exponents of the class groups in the said $p$-adic Lie extension. Our results generalize a previous result of Lei, where he established such an estimate under the assumption that $H\cong \Zp$.


\end{abstract}

\section{Introduction}

In this paper, $p$ will always denote an odd prime.  If $N$ is a finite abelian $p$-group, we shall write $e(N)= \log_p|N|$. Fix once and for all an algebraic closure $\bar{\Q}$ of the rational field $\Q$. Therefore, an algebraic (possibly infinite) extension of $\Q$ will mean an subfield of $\bar{\Q}$. A finite extension $F$ of $\Q$ is then said to be a number field. For every number field $F$, denote by $\Cl(F)$ the ideal class group of $F$. We can now state the following
celebrated asymptotic class number formula of Iwasawa \cite{Iw59} which is the main motivation behind this paper.

\begin{thm2}[Iwasawa]
Let $F_{\infty}$ be a $\Zp$-extension of a number field $F$. Denote $F_n$ to be the intermediate subfield of $F_{\infty}$ with index $|F_n:F|=p^n$. Then there exist $\mu, \lambda$ and $\nu$ (independent of $n$) such that
\[ e(\Cl(F_n)[p^{\infty}]) = \mu p^n+\lambda n +\nu\]
for $n\gg 0$.
\end{thm2}

Subsequently, this result was generalized by Cuoco and Monsky to the case of a $\Zp^d$-extension (cf. \cite[Theorem I]{CM}; also see \cite[Theorem 3.13]{Mon}). A common feature in the proofs of these formulas lies in the utilization of the structural theory of finitely generated modules over the commutative Iwasawa algebras.

In view of these results, the next natural direction of investigation is to consider the case of a noncommutative $p$-adic Lie extension $F_{\infty}/F$. Unfortunately, over a noncommutative Iwasawa algebra, one does not have a nice enough structural theory of modules to work with (see \cite{CFKSV, CSSalg}). One case where we do have such a structure theorem over noncommutative Iwasawa algebras, thanks to Howson \cite{Ho2} and Venjakob \cite{V02}, is when the module is finitely generated $p$-torsion. Building on this, Perbet \cite{Per} was able to prove certain results in this direction which we now describe. Let $F_{\infty}$ be a Galois extension of a number field $F$, whose Galois group $G=\Gal(F_{\infty}/F)$ is a compact pro-$p$ $p$-adic Lie group without $p$-torsion. Denote by $d$ the dimension of $G$. Let $\mathcal{M}$ be the maximal abelian unramified pro-$p$ extension of $F_{\infty}$. By maximality, $\mathcal{M}$ is a Galois extension of $F$. Write $\X=\Gal(\mathcal{M}/F_{\infty})$ and $\Y=\Gal(\mathcal{M}/F)$. The extension $1\lra \X\lra \Y\lra G\lra 1$ of groups induces a natural action of $G$ on $\X$ via conjugation by a lift in $\Y$ which in turn gives $\X$ a $\Zp\ps{G}$-module structure. It is well-known that $\X$ is a finitely generated $\Zp\ps{G}$-module (for instance, see \cite[Proposition 3.1]{Per}). By the work of Venjakob \cite{V02}, there is a notion of $\Zp\ps{G}$-rank of such modules. Furthermore, the works of Howson \cite{Ho, Ho2} and Venjakob \cite{V02} enable one to attach an Iwasawa $\mu_G$-invariant to $\X$. With these, Perbet was able to prove the following theorem which, for ease of exposition, we state in a slightly simplified form.

\begin{thm2}[Perbet]
Let $F_{\infty}$ be a $p$-adic-extension of a number field $F$ with Galois group being a uniform pro-$p$ group of dimension $d$. Denote by $F_n$ the intermediate subfield of $F_{\infty}$ with index $|F_n:F|=p^{dn}$. Then we have
\[ e(\Cl(F_n)[p^n]) =  \rank_{\Zp\ps{G}}(\X) n p^{dn} + \mu_G(\X) p^{dn} + O(np^{(d-1)n}).\]
\end{thm2}

We emphasis that Perbet's result is concerned with the growth of $\Cl(F_n)[p^n]$ rather than $\Cl(F_n)[p^{\infty}]$ as considered by Iwasawa and Cuoco-Monsky.

Recently, a great deal of research activities in noncommutative Iwasawa theory have been revolving around a $p$-adic Lie extension whose Galois group $G$ contains a closed normal subgroup $H$ such that $G/H\cong \Zp$ (for instance, see \cite{CFKSV, Ho, V03}). Following \cite{CFKSV}, we say that $\X$ satisfies the $\M_H(G)$-property if $\X_f: = \X/\X(p)$ is finitely generated over $\Zp\ps{H}$, where here $\X(p)$ is the $p$-primary submodule of $\X$. Note that in the event that $\X$ satisfies the $\M_H(G)$-property, it is then necessarily $\Zp\ps{G}$-torsion. Therefore, Perbet's theorem in this situation yields $e(\Cl(F_n)[p^n]) =  \mu_G(\X) p^{dn} + O(np^{(d-1)n})$.

The goal of this paper is to utilize structure theory of $p$-torsion modules to elucidate the error terms further. Our main result is as follows, where we succeed in at least obtaining an asymptotic upper bound.

\begin{thm2}[Theorem \ref{main arithmetic theorem}]
Let $F_{\infty}$ be a $p$-adic-extension of a number field $F$, where the Galois group $G=\Gal(F_{\infty}/F)$ is a uniform pro-$p$ group of dimension $d$ and contains a closed normal subgroup $H$ such that $G/H\cong \Zp$. Denote $F_n$ to be the intermediate subfield of $F_{\infty}$ with index $|F_n:F|=p^{dn}$. Suppose that $\X$ satisfies the $\M_H(G)$-property. Then we have
\[ e(\Cl(F_n)[p^n]) \leq \mu_G(\X) p^{dn} +  \rank_{\Zp\ps{H}}(\X_f) n p^{(d-1)n}+ O(p^{(d-1)n}).\]
\end{thm2}

In the case when $H\cong \Zp$, this was established by Lei \cite[Corollary 6.2]{Lei} under certain extra ramification conditions on the primes at $p$. Therefore, our result is a natural generalization of this previous result of Lei.

We should mention that only very recently, an asymptotic formula in the spirit of Iwasawa and Cuoco-Monsky was obtained for $\Zp^r\rtimes\Zp$-extension under the stronger assumption that $\X$ is finitely generated over $\Zp\ps{H}$ (see \cite{Lei, LiangLim}) plus some extra ramification conditions. In these works, the authors relies heavily on the structure theory of $\Zp\ps{H}$-modules, where in those situations one has $H\cong \Zp^r$. This unfortunately does not apply in our noncommutative situation considered in this paper. Despite this, it might be of interest to ask if we can still say anything in the noncommutative situation under the $\Zp\ps{H}$-finite generation assumption. This is the content of the next result, where we obtain an asymptotic upper bound with an error of $O(np^{(d-2)n})$.

\begin{thm2}[Theorem \ref{main arithmetic theorem2}]
Let $F_{\infty}$ be a $p$-adic-extension of a number field $F$, where the Galois group $G=\Gal(F_{\infty}/F)$ is a uniform pro-$p$ group of dimension $d$ and contains a closed normal subgroup $H$ such that $G/H\cong \Zp$. Denote by $F_n$ the intermediate subfield of $F_{\infty}$ with index $|F_n:F|=p^{dn}$. Suppose that $\X$ is finitely generated over $\Zp\ps{H}$. Then we have
\[ e(Cl(F_n)[p^n]) \leq \rank_{\Zp\ps{H}}(\X) n p^{(d-1)n}+ \mu_H(\X) p^{(d-1)n} +O(np^{(d-2)n}).\]
\end{thm2}

In the final section of the paper, we conclude with some examples. We mention that although some of the examples considered are $p$-adic Lie extensions with Galois group $\Zp\rtimes\Zp$, the results of Lei or Liang-Lim do not apply to these due to the extra ramification condition imposed in those results.

\begin{ack}
The author like to thank Antonio Lei for many insightful discussion on his paper \cite{Lei} and the subject on the asymptotic class number formulas in general. The author would also like to thank Dingli Liang for his interest and discussion on the subject of the paper. This research is supported by the
National Natural Science Foundation of China under Grant No.\ 11550110172 and Grant No.\ 11771164.
 \end{ack}

\section{Algebraic preliminaries}

As before, $p$ will denote an odd prime. Let $G$ be a compact
pro-$p$ $p$-adic Lie group with no $p$-torsion.
The completed group algebra of $G$ over $\Zp$ is then defined by
 \[ \Zp\ps{G} = \plim_U \Zp[G/U], \]
where $U$ runs over the open normal subgroups of $G$ and the inverse
limit is taken with respect to the canonical projection maps.  Under these said assumptions, it follows that
$\Zp\ps{G}$ is an Auslander regular ring (cf.
\cite[Theorem 3.26]{V02} or \cite[Theorem A.1]{LimFine}; for the definition of Auslander regular rings, see \cite[Definition 3.3]{V02}). Furthermore, the ring
$\Zp\ps{G}$ has no zero divisors (cf.\ \cite{Neu}), and therefore,
admits a skew field $Q(G)$ which is flat over $\Zp\ps{G}$ (see
\cite[Chapters 6 and 10]{GW} or \cite[Chapter 4, \S 9 and \S
10]{Lam}). For a finitely generated $\Zp\ps{G}$-module $M$,
its $\Zp\ps{G}$-rank is defined to be
$$ \rank_{\Zp\ps{G}}M  = \dim_{Q(G)} Q(G)\ot_{\Zp\ps{G}}M. $$
 We say that the $\Zp\ps{G}$-module $M$ is \textit{torsion} if
$\rank_{\Zp\ps{G}} M = 0$.  A finitely generated torsion $\Zp\ps{G}$-module $M$
is then said to be \textit{pseudo-null} if $\Ext^1_{\Zp\ps{G}}(M,
\Zp\ps{G}) =0$. Note that every subquotient of a
torsion $\Zp\ps{G}$-module (resp., pseudo-null $\Zp\ps{G}$-module)
is also torsion (resp., pseudo-null) (see \cite[Propositions 3.5(ii) and 3.6(ii)]{V02}).

Writing $\mathbb{F}_p$ for the finite field of order $p$, the completed group algebra of $G$ over $\mathbb{F}_p$ is given by
 \[ \mathbb{F}_p\ps{G} = \plim_U \mathbb{F}_p[G/U], \]
where $U$ runs over the open normal subgroups of $G$ and the inverse
limit is taken with respect to the canonical projection maps. Since we are assuming that
$G$ is pro-$p$ without $p$-torsion, it follows
that $\mathbb{F}_p\ps{G}$ is an Auslander
regular ring (cf. \cite[Theorem 3.30(ii)]{V02}) and has no zero divisors (cf. \cite[Theorem
C]{AB}).

We now define the notion of the Iwasawa $\mu_G$-invariant.
For a given finitely generated $\Zp\ps{G}$-module $M$, we write
$M(p)$ for the $\Zp\ps{G}$-submodule of $M$ which consists of
elements of $M$ that are annihilated by some power of $p$.
As seen in the previous paragraphs, the rings $\Zp\ps{G}$ and $\mathbb{F}_p\ps{G}$ are Auslander
regular and have no zero divisors. Therefore, we may apply \cite[Proposition
1.11]{Ho2} (see also \cite[Theorem 3.40]{V02}) to conclude that there is a
$\Zp\ps{G}$-homomorphism
\[ \varphi: M(p) \lra \bigoplus_{i=1}^s\Zp\ps{G}/p^{\al_i},\] whose
kernel and cokernel are pseudo-null $\Zp\ps{G}$-modules, and where
the integers $s$ and $\al_i$ are uniquely determined. The $\mu_G$-invariant of $M$ is then defined to be $\mu_G(M) = \displaystyle
\sum_{i=1}^s\al_i$.

Let $d$ denote the dimension of the group $G$. We shall once and for all fix an open normal uniform subgroup $G_0$ of $G$. Such a group exists by virtue of Lazard's theorem (see \cite[Corollary 8.34]{DSMS}). In the event that $G$ is already a uniform group, we shall take $G_0 = G$. We then write $G_n$ for the lower
$p$-series $P_{n+1}(G_0)$ which is defined recursively by $P_{1}(G_0) = G_0$, and
\[ P_{i+1}(G_0) = \overline{P_{i}(G_0)^{p}[P_{i}(G_0),G_0]}, ~\mbox{for}~ i\geq 1. \]
It follows from \cite[Thm.\
3.6]{DSMS} that $G^{p^i} =
P_{i+1}(G)$ and that we have an equality $|G_0:P_2(G_0)| = |P_i(G_0):
P_{i+1}(G_0)|$ for every $i\geq 1$ (cf. \cite[Definition 4.1]{DSMS}). It is not difficult to
verify that $|G:G_n| = [G:G_0]p^{dn}$, where $d= \dim G$. We now record the following lemma whose proof is left to the readers as an exercise.

\bl \label{lemma rank mu}
 Let $M$ be a finitely generated $\Zp\ps{G}$-module. Then $M$ is finitely generated over $\Zp\ps{G_0}$ with
 \[\rank_{\Zp\ps{G_0}}(M) =[G:G_0]\rank_{\Zp\ps{G}}(M)\quad \mbox{and}\quad
  \mu_{G_0}(M) = [G:G_0]\mu_G(M).\]
\el

From now on, we further suppose that the group $G$ contains a closed normal subgroup $H$ with the property that $\Ga:=G/H\cong \Zp$. Since $G_0$ is an open uniform subgroup of $G$, $H_0:=H\cap G_0$ is also an open uniform subgroup of $H$. Write $H_n$ for the lower $p$-series $P_{n+1}(H_0)$ of $H_0$. Set $\Ga_0 = G_0/H_0$ and write $\Ga_n= \Ga_0^{p^n}$.

\bl
For every $n\geq 1$, we have $H_n = H\cap G_n$ and $G_n/H_n\cong\Ga_n$.
\el

\bpf
It clearly suffices to prove the lemma for the case $G=G_0$ is uniform. Then since $H$ and $G$ are uniform, we have  $H_n=H^{p^n}$ and $G_n = G^{p^n}$ (cf. \cite[Theorem 3.6]{DSMS}). Clearly, we have $H^{p^n}\subseteq H\cap G^{p^n}$. Conversely, let $h\in H\cap G^{p^n}$. Then there exists $g\in G$ such that $h= g^{p^n}$ which in turn implies that the coset $gH$ is a torsion element in $G/H$. But since $G/H\cong\Zp$ has no $p$-torsion, we have $g\in H$, and hence $h\in H^{p^n}$. This proves the first equality. For the second equality, we simply observe that \[G_n/H_n = G^{p^n}/H^{p^n}\cong G^{p^n}H/H= (G/H)^{p^n}\cong \Ga^{p^n} = \Ga_n.\]
\epf

Following \cite{CFKSV}, we say that a finitely generated $\Zp\ps{G}$-module $M$ satisfies the $\M_H(G)$-property if $M_f := M/M(p)$ is finitely generated over $\Zp\ps{H}$.
For a finite $\Zp$-module $N$, we write $e(N)$ for the $p$-exponent of $N$, i.e., $|N| = p^{e(N)}$. We can now state the main algebraic results of this section.

\bp \label{main alg theorem}
 Let $G$ be a compact pro-$p$ $p$-adic Lie group without $p$-torsion. Suppose that $G$ contains a closed normal subgroup $H$ with the property that $\Ga:=G/H\cong \Zp$. Let $M$ be a finitely generated $\Zp\ps{G}$-module which satisfies the $\M_H(G)$-property. Then we have
 \[ e(M_{G_n}/p^n) \leq \mu_G(M)p^{dn} +\rank_{\Zp\ps{H}}(M_f)np^{(d-1)n} + O(p^{(d-1)n}).\]
\ep

\bpf
 In view of Lemma \ref{lemma rank mu}, it suffices to prove the theorem under the assumption that $G$ (and hence $H$) is a uniform pro-$p$ group with $G_n\cap H= H_n$ which we will do. Note that under this assumption, we have $\Ga_n=G_n/H_n$. Now since
the ring $\Zp\ps{G}$ is Noetherian, the submodule $M(p)$ is certainly finitely
generated over $\Zp\ps{G}$. Therefore, one can find an integer
$t$ such that $p^t$ annihilates $M(p)$. Let $n\geq t$. Consider the following commutative diagram
\[  \entrymodifiers={!! <0pt, .8ex>+} \SelectTips{eu}{}\xymatrix{
    0 \ar[r]^{} & M(p) \ar[d]_{p^n} \ar[r] &
    M \ar[d]_{p^n}
    \ar[r]^{} & M_f \ar[d]_{p^n} \ar[r] &0\\
    0 \ar[r]^{} & M(p) \ar[r]^{}
    & M \ar[r] & M_f \ar[r] &0 } \]
with exact rows. Note that the rightmost vertical map is injective, and by our choice of $n$, the leftmost vertical map is zero. The snake lemma therefore gives a short exact sequence
\[ 0\lra M(p) \lra M/p^n \lra M_f/p^n\lra 0. \]
As $G_n$ is an open subgroup of $G$, the above short exact sequence is also a short exact sequence of finitely generated $\Zp\ps{G_n}$-module. Upon taking $G_n$-invariant, we obtain an exact sequence
\[ M(p)_{G_n} \lra (M/p^n)_{G_n} \lra (M_f/p^n)_{G_n}\lra 0 \]
of finitely generated $\Zp$-modules (cf. \cite[Lemma 3.2.3]{LS}). On the other hand, since every module appearing in this exact sequence is annihilated by $p^n$, the exact sequence is a sequence of finite $\Zp$-modules. Hence we have an inequality
\[ e((M/p^n)_{G_n})\leq e(M(p)_{G_n})+ e((M_f/p^n)_{G_n}).\]
By \cite[Theorem 2.5.1]{LimCMu}, we have
\[  e(M(p)_{G_n})=\mu_G(M)p^{dn} + O(p^{(d-1)n}).\]
It therefore remains to estimate $e((M_f/p^n)_{G_n})$. Since $M_f$ is finitely generated over $\Zp\ps{H}$, we may apply \cite[Theorem 2.1(ii)]{Per} to obtain
\[  e((M_f/p^n)_{H_n})=\rank_{\Zp\ps{H}}(M_f)np^{(d-1)n} + O(np^{(d-2)n}),\]
noting that $H$ has dimension $d-1$ and $\mu_H(M_f)=0$ since $M_f(p)=0$. By the $\Zp\ps{H}$-finite generation of $M_f$, we have that $(M_f/p^n)_{H_n}$ is finite by a similar argument as above. It then follows from this that we have an inequality
\[e((M_f/p^n)_{G_n}) =  e(((M_f/p^n)_{H_n})_{\Ga_n}) \leq e((M_f/p^n)_{H_n}).\]
Combining all these estimates, we have the required bound of the proposition.
\epf

The next result is a variant of the preceding result which gives an estimate under a stronger hypothesis on the structure of our module $M$.

\bp \label{main alg theorem2}
 Let $G$ be a compact pro-$p$ $p$-adic Lie group without $p$-torsion. Suppose that $G$ contains a closed normal subgroup $H$ with the property that $\Ga:=G/H\cong \Zp$. Let $M$ be a $\Zp\ps{G}$-module which is finitely generated over $\Zp\ps{H}$. Then we have
 \[ e(M_{G_n}/p^n) \leq \rank_{\Zp\ps{H}}(M)np^{(d-1)n} + \mu_H(M)p^{(d-1)n} + O(np^{(d-2)n}).\]
\ep

\bpf
By a similar argument to that in Proposition \ref{main alg theorem}, we have
\[ e((M/p^n)_{G_n})\leq e(M(p)_{G_n})+ e((M_f/p^n)_{G_n})\]
and
\[e((M_f/p^n)_{G_n}) \leq \rank_{\Zp\ps{H}}(M)np^{(d-1)n}+ O(np^{(d-2)n}).\] On the other hand, we have
\[ e(M(p)_{G_n})=e(M(p)_{H_n})_{\Ga_n}\leq e(M(p)_{H_n}) = \mu_H(M)p^{(d-1)n} + O(p^{(d-2)n}),\]
where the last equality follows from an application of \cite[Theorem 2.5.1]{LimCMu}.
Combining these two estimates, we obtain the required bound of the theorem.
\epf

We end the section with the another useful estimate.

\bl \label{Z/p estimate}
 Let $G$ be a compact pro-$p$ $p$-adic Lie group without $p$-torsion. Write $d=\dim G$. Then we have
 \[ e(H_1(G_n, \Z/p^n))\leq dn\quad \mbox{and}\quad e(H_2(G_n, \Z/p^n))\leq  {d\choose 2}n.\]
\el

\bpf
 By mathematical induction and the long cohomology sequence of the short exact sequence
 \[ 0\lra \Z/p^{i}\lra \Z/p^{i+1}\lra \Z/p\lra 0, \]
 we have
  \[ e(H_1(G_n, \Z/p^n))\leq n  e(H_1(G_n, \Z/p)))\]
  and
  \[ e(H_2(G_n, \Z/p^n))\leq n  e(H_2(G_n, \Z/p))).\]
  Since $G_n$ is a uniform group of dimension $d$, $e(H_1(G_n, \Z/p)))=d$ and $e(H_2(G_n, \Z/p)))={d\choose 2}$ (cf. \cite[Theorem 4.35]{DSMS}). This proves the inequalities of the proposition.
\epf

\section{Arithmetic setup and the main theorem}

We turn to arithmetic. Fix once and for all an algebraic closure $\bar{\Q}$ of $\Q$. Therefore, an algebraic (possibly infinite) extension of $\Q$ will mean an subfield of $\bar{\Q}$. A finite extension $F$ of $\Q$ is called a number field, and we fix one such $F$ as our base field. Let $F_{\infty}$ be a $p$-adic Lie extension of $F$ with Galois group $G$. Suppose that $G$ is pro-$p$ torsionfree and contains a closed normal subgroup $H$ such that $\Ga:=G/H\cong\Zp$. We shall further assume that our extension $F_{\infty}/F$ satisfies the following ramification property.

\medskip
\noindent \textbf{(Ram$_{S}$):}  $F_{\infty}$ is unramified outside a finite set of primes of $F$.

\medskip
Let $\mathcal{M}$ be the maximal abelian unramified pro-$p$ extension of $F_{\infty}$. By maximality, $\mathcal{M}$ is also Galois over $F$. Write $\X=\Gal(\mathcal{M}/F_{\infty})$ and $\Y=\Gal(\mathcal{M}/F)$. Clearly, $\Y/\X\cong G$. There is a natural action of $G$ on $\X$ defined as follows: for $x\in \X$ and $g\in G$, $x^g:=\widetilde{g}x\widetilde{g}^{-1}$, where $\widetilde{g}\in \Y$ is a lift of $g$. It is well-known that $\X$ is a finitely generated $\Zp\ps{G}$-module (cf. \cite[Proposition 3.1]{Per}).

As in Section 2, we shall fix a normal uniform subgroup $G_0$ of $G=\Gal(F_{\infty}/F)$ and write $G_n$ for the lower $p$-series of $G_0$. The corresponding fixed field of $G_n$ is then denoted to be $F_n$. We can now state the following theorem which generalizes \cite[Corollary 6.2]{Lei}.

\bt \label{main arithmetic theorem}
Let $F_{\infty}$ be a $p$-adic Lie extension of a number field $F$ with Galois group $G$. Suppose that the following conditions are all valid.
\begin{enumerate}
\item[$(a)$] $G$ is a pro-$p$ torsionfree $p$-adic Lie group of dimension $d\geq 2$.
\item[$(b)$] $G$ contains a closed normal subgroup $H$ with $G/H\cong \Zp$.
\item[$(c)$] $\X$ satisfies the $\M_H(G)$-property.
\item[$(d)$] \textbf{(Ram$_{S}$)} is valid.
\end{enumerate}
Then \[e\big(\Cl(F_n)[p^n]\big)\leq \mu_G(\X)p^{dn} +\rank_{\Zp\ps{H}}(\X_f)np^{(d-1)n} + O(p^{(d-1)n}).\]
\et

\bpf
From the spectral sequence
\[ H_r(G_n,H_s(\X,\Z/p^n)) \Longrightarrow H_{r+s}(\Y_n, \Z/p^n), \]
we have
\[ H_2(G_n,\Z/p^n)\lra (\X/p^n)_{G_n}\lra \Y_n^{ab}/p^n\lra H_1(G_n,\Z/p^n)\lra 0. \]
By virtue of Proposition \ref{main alg theorem} and Lemma \ref{Z/p estimate}, we have
\[e\big(\Y_n^{ab}/p^n\big)\leq \mu_G(\X)p^{dn} +\rank_{\Zp\ps{H}}(\X_f)np^{(d-1)n} + O(p^{(d-1)n}).\]
On the other hand, class field theory gives us a short exact sequence
\[ 0\lra \bar{C}_n \lra \Y_n^{ab} \lra \Cl(F_n)[p^{\infty}]\lra 0  \]
which in turn induces the following exact sequence
\[  \bar{C}_n/p^n \lra \Y_n^{ab}/p^n \lra \Cl(F_n)[p^{\infty}]/p^n\lra 0.  \]
It then follows from our estimate for $\Y_n^{ab}/p^n$ and this exact sequence that
\[ e\big(\Cl(F_n)[p^{\infty}]/p^n\big) \leq e\big(\Y_n^{ab}/p^n\big)\leq \mu_G(\X)p^{dn} +\rank_{\Zp\ps{H}}(\X_f)np^{(d-1)n} + O(p^{(d-1)n}).\]
But since $\Cl_n(F)$ is finite, we have $e(\Cl(F_n)[p^{n}])=e\big(\Cl(F_n)[p^{\infty}]/p^n\big)$ and this proves the theorem.
\epf

\medskip
The next theorem is a variant of the previous, where we can elucidate the error term under a stronger assumption on the structure of $\X$.

\bt \label{main arithmetic theorem2}
Let $F_{\infty}$ be a $p$-adic Lie extension of a number field $F$ with Galois group $G$. Suppose that the following conditions are all valid.
\begin{enumerate}
\item[$(a)$] $G$ is a pro-$p$ torsionfree $p$-adic Lie group of dimension $d\geq 2$.
\item[$(b)$] $G$ contains a closed normal subgroup $H$ with $G/H\cong \Zp$.
\item[$(c)$] $\X$ is finitely generated over $\Zp\ps{H}$.
\item[$(d)$] The assumption \textbf{(Ram$_{S}$)} is valid.
\end{enumerate}
Then \[e\big(\Cl(F_n)[p^n]\big)\leq \rank_{\Zp\ps{H}}(\X)np^{(d-1)n} + \mu_H(\X)p^{(d-1)n}+O(np^{(d-2)n}).\]
\et

\bpf
This has a similar proof to that of Theorem \ref{main arithmetic theorem}, where we made use of Proposition \ref{main alg theorem2} in place of Proposition \ref{main alg theorem}.
\epf

\section{Examples}

We end the paper discussing some examples to illustrate the results of this paper.

Let $F=\Q(\mu_p)$ and $K=\Q(\mu_{p^{\infty}})$. Denote by $\mathcal{M}_K$ the maximal abelian unramified pro-$p$ extension of $K$. Write $\X_K$ for $\Gal(\mathcal{M}_K/K)$. Let $\la_K$ denote the $\Zp$-rank of $\X_K^-$, where here $\X_K^-$ is the minus one eigenspace under complex conjugation. As seen in \cite[Example 5.1]{HS}, there exist at least $\la_K$ $\Zp$-extensions of $K$ which is unramified outside $p$. Let $F_{\infty}$ be one of such $\Zp$-extension of $K$. Then $\Gal(F_{\infty}/F)\cong \Zp\rtimes\Zp$. Denote by $\X=\X_{F_{\infty}}$ the maximal abelian unramified pro-$p$ extension of $F_{\infty}$. This is a finitely generated $\Zp\ps{H}$-module of $\Zp\ps{H}$-rank $\geq \la_K-1$, where $H=\Gal(F_{\infty}/K)$ (cf. \cite[Theorem 4.1 and Example 5.1]{HS}). Therefore, we may apply Theorem \ref{main arithmetic theorem2} to obtain an asymptotic upper bound for the class groups in this tower. We mention that the work of Lei \cite{Lei} and Liang-Lim \cite{LiangLim} do not apply here, since the results there only applied to $p$-adic extensions which are totally ramified at the prime $p$ and the extension considered clearly do not satisfy this condition.

Similarly, we can also apply Theorem \ref{main arithmetic theorem2} to
\cite[Examples 5.2 and 5.3]{HS}.



\end{document}